\documentclass[12pt,reqno]{amsart}

\usepackage{amssymb}

\usepackage{eucal}

\input{diagrams}

\font\sss=cmss8

\def\BQ{{\mathbb Q}}

\def\BZ{{\mathbb Z}}

\def\sC{\mbox{\sf C}}
\def\sD{\mbox{\sf D}}

\def\ast{{\textstyle *}}
\def\Aut{\operatorname{Aut}}

\def\C{\operatorname{C}}

\def\D{\sD}
\def\Dc{\sD^{\operatorname{c}}}

\def\Df{\sD^{\operatorname{f}}}
\def\dim{\operatorname{dim}}
\def\Dsmall{\mbox{\sss D}}
\def\dual{\operatorname{D}}
\def\End{\operatorname{End}}
\def\Ext{\operatorname{Ext}}

\def\H{\operatorname{H}}

\def\Hom{\operatorname{Hom}}

\def\inf{\operatorname{inf}}

\def\LTensor{\stackrel{\operatorname{L}}{\otimes}}
\def\opp{\operatorname{op}}

\def\RHom{\operatorname{RHom}}

\def\Tor{\operatorname{Tor}}

\numberwithin{equation}{part}

\newtheorem{Lemma}{Lemma}[section]
\newtheorem{Theorem}[Lemma]{Theorem}

\newtheorem{Corollary}[Lemma]{Corollary}

\theoremstyle{definition}

\newtheorem{Setup}[Lemma]{Setup}

\def\AR{Aus\-lan\-der-Rei\-ten}
\def\ARq{\AR\ qui\-ver}
\def\ARt{\AR\ tri\-ang\-le}
\def\ARtr{\AR\ trans\-la\-tion}
\def\compact{small}

\def\DGm{DG mo\-du\-le}

\def\DGRlm{DG left-$R$-mo\-du\-le}

\def\DGRrm{DG right-$R$-mo\-du\-le}

\def\DGRbm{DG left/right-$R$-mo\-du\-le}

\def\DGRolm{DG left-$R^{\opp}$-mo\-du\-le}

\def\grlm{graded left-mo\-du\-le}

\def\DGCXklm{DG left-$\CXk$-module}

\def\Space{X}  

\def\CXk{\C^{\ast}(\Space;k)}

\def\HXk{\H^{\ast}(\Space;k)}

\def\DcR{\Dc(R)}
\def\DcRo{\Dc(R^{\opp})}
\def\DfR{\Df(R)}
\def\DfRo{\Df(R^{\opp})}
\def\DcXk{\Dc(\CXk)}

\def\Tree{T}
\def\Group{\Pi}
\def\Component{C}

\def\AddFunct{\beta}
\def\AddFunctTree{f}

\begin{document}

\title[Auslander-Reiten quiver]
{The Auslander-Reiten quiver of a Poincar\'e duality space}

\author{Peter J\o rgensen}
\address{Danish National Library of Science and Medicine, N\o rre
All\'e 49, 2200 K\o \-ben\-havn N, DK--Denmark}
\email{pej@dnlb.dk, www.geocities.com/popjoerg}




\begin{abstract} 

In a previous paper, \ARt s and quivers were introduced into algebraic
topology.  This paper shows that over a Poincar\'e duality space,
each component of the \ARq\ is isomorphic to $\BZ A_{\infty}$.

\end{abstract}

\maketitle

\setcounter{section}{-1}
\section{Introduction}
\label{sec:introduction}

In \cite{PJARtq}, the concepts of \ARt s and \ARq s from the
representation theory of Artin algebras were introduced into algebraic
topology.  

The main theorem was that \ARt s exist precisely over Poincar\'e
duality spaces.  On the other hand, the only concrete \ARq s computed
in \cite{PJARtq} were those over spheres.  An ad hoc computation
showed that the \ARq\ over $S^d$, the $d$-dimensional sphere, consists
of $d-1$ components, each isomorphic to $\BZ A_{\infty}$.

The purpose of this paper is to show that this result was no accident.
Namely, let $k$ be a field of characteristic zero, and let $\Space$ be
a simply connected topological space which has Poincar\'e duality of
dimension $d \geq 2$ over $k$, that is, satisfies
$\Hom_k(\H^{\ast}(\Space;k),k) \cong \H^{d-\ast}(\Space;k)$, where
$\HXk$ is singular cohomology of $X$ with coefficients in $k$.  The
main result here is the following.

\begin{Theorem}
\label{thm:topology}
Let $\Component$ be a component of the \ARq\ of the category $\DcXk$.
Then $\Component$ is isomorphic to $\BZ A_{\infty}$.
\end{Theorem}

Here $\CXk$ is the singular cochain Differential Graded Algebra of
$\Space$ with coefficients in $k$, and if $\D(\CXk)$ is the
derived category of Differential Graded left-modules over $\CXk$, then
$\DcXk$ denotes the full subcategory of \compact\ objects, that
is, objects $M$ for which the functor $\Hom(M,-)$ preserves set
indexed coproducts.

The ingredients in the proof are standard.  First the Riedtmann
structure theorem is invoked to show that the component $\Component$
has the form $\BZ \Tree/\Group$ for a directed tree $\Tree$ and an
admissible group of automorphisms $\Group \subseteq \Aut(\BZ \Tree)$
(section \ref{sec:Riedtmann}).  Then a certain function $\AddFunct$ is
constructed on $\Component$ (section \ref{sec:AddFunct}).  This
induces an unbounded additive function $\AddFunctTree$ on $\Tree$
forcing its underlying graph to be $A_{\infty}$, and finally an
elementary argument shows that $\Group$ acts trivially (section
\ref{sec:shape}).  Hence $\Component$ is $\BZ A_{\infty}$.

In fact, I will prove a more general result than theorem
\ref{thm:topology}.  Namely, I will not restrict myself to 
Differential Graded Algebras (DGAs) of the form $\CXk$, but will work
over an abstract DGA denoted $R$.

\begin{Setup}
\label{set:blanket}
Throughout, $R$ denotes a DGA over the field $k$ which satisfies
the following.
\begin{enumerate}

  \item  $R$ is a cochain DGA, that is, $R^i = 0$ for $i < 0$.

  \item  $R^0 = k$.

  \item  $R^1 = 0$.

  \item  $\dim_k R < \infty$. \hfill $\Box$

\end{enumerate}
\end{Setup}

This setup is the same as in \cite[sec.\ 3]{PJARtq}, so I can use the
results about $R$ proved in \cite{PJARtq}.

For brief introductions to \ARt s, \ARq s, DGAs, Differential Graded
modules (\DGm s), and the way they interact, I refer the reader to
\cite{PJARtq}.  My notation is mostly standard and identical to
the notation of \cite{PJARtq}; however, I do want to recapitulate a
few ubiquitous items.

By $\D(R)$ is denoted the derived category of \DGRlm s, and by $\DcR$
the full subcategory of \compact\ objects.

By $R^{\opp}$ is denoted the opposite DGA of $R$ with multiplication
\[
  r \stackrel{\opp}{\cdot} s = (-1)^{|r||s|}sr.  
\]
\DGRolm s are identified with \DGRrm s, and so $\D(R^{\opp})$ is
identified with the derived category of \DGRrm s, and $\DcRo$ with its
full subcategory of \compact\ objects.

The \DGRbm s
\[
  k = R/R^{\geq 1} \; \; \mbox{and} \; \; \dual\!R = \Hom_k(R,k) 
\]
are used frequently throughout.

When I wish to emphasize that I am viewing some $M$ as either a
\DGRlm\ or a \DGRrm, I do so with subscripts, writing either ${}_{R}M$
or $M_R$.

\section{Applying the Riedtmann structure theorem}
\label{sec:Riedtmann}

This section will apply the Riedtmann structure theorem to show under
some conditions that a component $\Component$ of the \ARq\ of $\DcR$
has the form $\BZ \Tree/\Group$ for a directed tree $\Tree$ and an
admissible group of automorphisms $\Group \subseteq \Aut(\BZ \Tree)$.

Let me start by recalling some facts from \cite{PJARtq}.  The
triangulated category $\DcR$ does not necessarily have \ARt s, but if
it does, then by \cite[prop.\ 4.4(ii)]{PJARtq} its \ARtr\ is
\begin{equation}
\label{equ:tau}
  \tau(-) = \Sigma^{-1}(\dual\!R \LTensor_R -).
\end{equation}
Here $\Sigma$ denotes suspension of \DGm s.

When $\DcR$ has \ARt s, $\tau$ induces a map from the \ARq\ of $\DcR$
to itself which is also referred to as the \ARtr\ and denoted $\tau$.

The \ARq\ of $\DcR$ is then a stable translation quiver with
translation $\tau$; see \cite[sec.\ 2]{PJARtq}.  This is even true in
the strong sense that $\tau$ is an automorphism of the \ARq, under the
extra assumption that $\DcRo$ also has \ARt s.  This last fact was not
mentioned in \cite{PJARtq}, but is proved in the following lemma.

\begin{Lemma}
\label{lem:tau_automorphism}
Suppose that $\DcR$ and $\DcRo$ have \ARt s.
\begin{enumerate}

  \item  The functor $\Sigma^{-1}(\dual\!R \LTensor_R -)$ is an
         auto-equivalence of $\DcR$, with quasi-inverse
         $\Sigma(P \LTensor_R -)$ for a suitable \DGRbm\ $P$.

  \item  The \ARtr\ $\tau$ is an automorphism of the \ARq\ of $\DcR$.

  \item  For each integer $p \not= 0$, the map $\tau^p$ is without fixed
         points in the \ARq\ of $\DcR$. 

\end{enumerate}
\end{Lemma}

\begin{proof}
(i)  Both $\DcR$ and $\DcRo$ have \ARt s, so it follows from
\cite[cor.\ 5.2 and thm.\ 5.1 and its proof]{PJARtq} that there are
isomorphisms ${}_{R}(\dual\!R) \cong {}_{R}(\Sigma^d R)$ in $\DcR$ and
$(\dual\!R)_R \cong (\Sigma^d R)_R$ in $\DcRo$.  (Note that I do not
know ${}_{R}(\dual\!R)_R \cong {}_{R}(\Sigma^d R)_R$!)

This makes it easy to check that with $P$ equal to
$\RHom_R(\dual\!R,R)$, the functors $\Sigma^{-1}(\dual\!R \LTensor_R
-)$ and $\Sigma(P \LTensor_R -)$ are quasi-inverse endofunctors on
$\DcR$.

\medskip
\noindent
(ii)  This clearly follows from (i).

\medskip
\noindent
(iii)  Let $M$ be an indecomposable object of $\DcR$ with vertex $[M]$
in the \ARq\ of $\DcR$.  If $\tau^p$ had the fixed point $[M]$, then
$p$ applications of equation \eqref{equ:tau} would give
\[
  M \cong \Sigma^{-p}((\dual\!R) \LTensor_R \cdots 
                      \LTensor_R (\dual\!R) \LTensor_R M).
\]
But this is impossible, as one proves for instance by checking 
\begin{align*}
  & \inf \{\, i \,|\, \H^i(M) \not= 0 \,\} \\
  & \;\;\;\;\;\;\;\;\;\; \not= 
    \inf \{\, i \,|\, \H^i(\Sigma^{-p}((\dual\!R) 
                           \LTensor_R \cdots \LTensor_R (\dual\!R)
                           \LTensor_R M)) \not= 0 \,\},
\end{align*}
for which the general formula 
\begin{align*}
  & \inf \{\, i \,|\, \H^i(N \LTensor_R M) \not= 0 \,\} \\
  & \;\;\;\;\;\;\;\;\;\; = 
    \inf \{\, i \,|\, \H^i(N) \not= 0 \,\}
    + \inf \{\, i \,|\, \H^i(M) \not= 0 \,\}
\end{align*}
is handy.
\end{proof}

A loop in a quiver is an arrow that starts and ends at the same
vertex.

\begin{Lemma}
\label{lem:no_loops}
Suppose that $\DcR$ has \ARt s.  In this case, the \ARq\ of $\DcR$ has no
loops. 
\end{Lemma}

\begin{proof}
Suppose that there were an indecomposable object $M$ of $\DcR$ whose
vertex $[M]$ in the \ARq\ of $\DcR$ had a loop $[M]
\longrightarrow [M]$.  The arrow $[M] \longrightarrow [M]$ would mean
that there existed an irreducible morphism $M
\stackrel{\mu}{\longrightarrow} M$ in $\DcR$.  Such a morphism would
be non-invertible, so it would be in the radical of the local artinian
ring $\End_{\Dsmall^{\operatorname{c}}(R)}(M)$ (cf.\
\cite[lem.\ 3.6]{PJARtq}).  Hence some power $\mu^n$ would be zero.

However, it is not hard to mimick the proof of \cite[lem.\
VII.2.5]{ARSbook} to show that if $M$ is indecomposable and $M
\stackrel{\mu}{\longrightarrow} M$ is irreducible, and the composition
$M \stackrel{\mu}{\longrightarrow} \cdots
\stackrel{\mu}{\longrightarrow} M$ is zero, then $[M] = \tau[M]$
holds, and this is false by lemma \ref{lem:tau_automorphism}(iii).
\end{proof}

In the following lemma and the rest of the paper, $\Component$ is a
component of the \ARq\ viewed as a stable translation quiver, so
$\Component$ is closed under the \ARtr\ $\tau$ and its inverse.  This
implies that $\Component$ is a connected stable translation quiver
with translation the restriction of $\tau$.

For the lemma's terminology of directed trees, stable translation
quivers of the form $\BZ \Tree$, and admissible groups of
automorphisms, see \cite[sec.\ 4.15]{Bensonbook}.

\begin{Lemma}
\label{lem:Riedtmann}
Suppose that $\DcR$ and $\DcRo$ have \ARt s.  Let $\Component$ be a
component of the \ARq\ of $\DcR$.

Then there exists a directed tree $\Tree$ and an admissible group of
automorphisms $\Group \subseteq \Aut(\BZ \Tree)$ so that $\Component$
is isomorphic to $\BZ \Tree/\Group$ as a stable translation
quiver.
\end{Lemma}

\begin{proof}
Lemma \ref{lem:tau_automorphism}(ii) implies that the restriction of
$\tau$ to $\Component$ is an automorphism of $\Component$, and lemma
\ref{lem:no_loops} implies that $\Component$ has no loops.  Also, by
construction, the \ARq\ has no multiple arrows, so neither has
$\Component$.

Hence $\Component$ satisfies the conditions of the Riedtmann structure
theorem, \cite[thm.\ 4.15.6]{Bensonbook}, and the conclusion of the
lemma follows.
\end{proof}

\section{Labelling the \ARq}
\label{sec:labelling}

This section constructs a labelling of the \ARq\ of $\DcR$, that is, a
pair of positive integers $(a_{m \rightarrow n},b_{m \rightarrow n})$
for each arrow $m \longrightarrow n$ in the quiver.  The construction
is classical, works in good cases such as when $\DcR$ and $\DcRo$ both
have \ARt s, and goes as follows.

Let $M$ and $N$ be indecomposable objects of $\DcR$ for which there is
an arrow $[M] \longrightarrow [N]$ in the \ARq, and let
\[
  M \longrightarrow X \longrightarrow \tau^{-1}M \longrightarrow
  \; \; \mbox{and} \; \;
  \tau N \longrightarrow Y \longrightarrow N \longrightarrow
\]
be \ARt s in $\DcR$.  Here $\tau^{-1}M$ makes sense because of lemma
\ref{lem:tau_automorphism}.  Note that the triangles are determined up
to isomorphism by \cite[prop.\ 3.5(i)]{HapDerCat}.

The arrow $[M] \longrightarrow [N]$ means that there exists an
irreducible morphism $M \longrightarrow N$ in $\DcR$, and by
\cite[prop.\ 3.5]{HapDerCat} this means that $M$ is a direct summand
of $Y$ and that $N$ is a direct summand of $X$.  Let $a$ be the
multiplicity of $M$ as a direct summand of $Y$ and let $b$ be the
multiplicity of $N$ as a direct summand of $X$.  These numbers are
well-defined because $\DcR$ is a Krull-Schmidt category by \cite[par.\
3.1]{HapDerCat} and \cite[lem.\ 3.6]{PJARtq}.  Now define the
labelling by equipping $[M] \longrightarrow [N]$ with the pair of
positive integers $(a,b)$.

\begin{Lemma}
\label{lem:symmetry}
Suppose that $\DcR$ and $\DcRo$ have \ARt s.  The above labelling of
the \ARq\ of $\DcR$ satisfies the following.
\begin{enumerate}

  \item  $a_{\tau[M] \rightarrow [N]} = b_{[N] \rightarrow [M]}$.

  \item  $b_{\tau[M] \rightarrow [N]} = a_{[N] \rightarrow [M]}$.

  \item  $(a_{\tau^p[M] \rightarrow \tau^p[N]},
           b_{\tau^p[M] \rightarrow \tau^p[N]}) =
          (a_{[M] \rightarrow [N]},b_{[M] \rightarrow [N]})$ for each 
         integer $p$.

\end{enumerate}
\end{Lemma}

\begin{proof}
(i)  Let 
\begin{equation}
\label{equ:AR_triangle}
  \tau N \longrightarrow Y \longrightarrow N \longrightarrow
\end{equation}
be an \ARt.  Then $a_{\tau[M] \rightarrow [N]}$ is defined as the
multiplicity of $\tau M$ as a direct summand of $Y$.  

Now, $\tau$ is given by the formula $\Sigma^{-1}(\dual\!R \LTensor_R
-)$ which also defines an auto-equivalence of $\DcR$ by lemma
\ref{lem:tau_automorphism}(i).  A quasi-inverse auto-equivalence is 
given by $\Sigma(P \LTensor_R -)$, and applying this to the \ARt\
\eqref{equ:AR_triangle} gives an \ARt\ 
\[
  N \longrightarrow \Sigma(P \LTensor_R Y) \longrightarrow \tau^{-1}N
  \longrightarrow.
\]
Then $b_{[N] \rightarrow [M]}$ is defined as the multiplicity of $M$
as a direct summand of $\Sigma(P \LTensor_R Y)$.  

Applying $\Sigma^{-1}(\dual\!R \LTensor_R -)$, this multiplicity
equals the multiplicity of $\Sigma^{-1}(\dual\!R \LTensor_R M)$ as a
direct summand of $\Sigma^{-1}(\dual\!R \LTensor_R \Sigma(P \LTensor_R
Y))$, that is, the multiplicity of $\tau M$ as a direct summand of
$Y$.  And this again equals $a_{\tau[M] \rightarrow [N]}$ by the
above.

\medskip
\noindent
(ii) and (iii) are proved by similar means.
\end{proof}

\section{The function $\AddFunct$}
\label{sec:AddFunct}

Let
\[
  \AddFunct(M) = \dim_k \Ext_R(M,k) 
  = \dim_k \H(\RHom_R(M,k))
\]
for $M$ in $\D(R)$.  Observe that $\AddFunct$ is constant on each
isomorphism class in $\DcR$, so induces a well-defined function on
the vertices $[M]$ of the \ARq\ of $\DcR$.  I also denote the induced
function by $\AddFunct$, so $\AddFunct([M]) = \AddFunct(M)$.

The purpose of this section is to study $\AddFunct$.

\begin{Lemma}
\label{lem:AddFunct}
Suppose that $\DcR$ and $\DcRo$ have \ARt s and that the \DGRlm\
${}_{R}k$ is not in $\DcR$.  

Let $N$ be an indecomposable object of $\DcR$ with vertex $[N]$ in the
\ARq\ of $\DcR$.
\begin{enumerate}

  \item  If $\tau N \longrightarrow Y \longrightarrow N 
\longrightarrow$ 
is an \ARt\ in $\DcR$, then 
$\AddFunct(\tau N) + \AddFunct(N) - \AddFunct(Y) = 0$.

  \item  $\AddFunct(\tau[N]) + \AddFunct([N]) 
- \sum_{[M] \rightarrow [N]}a_{[M] \rightarrow [N]}\AddFunct([M]) = 0$,
where the sum is over all arrows in the \ARq\ of $\DcR$ which end in
$[N]$. 

  \item  $\AddFunct(\tau[N]) = \AddFunct([N])$.

\end{enumerate}
\end{Lemma}

\begin{proof}
(i)  The \ARt\ $\tau N \longrightarrow Y \longrightarrow N
\longrightarrow$ gives a long exact sequence consisting of pieces
$\Ext_R^i(N,k) \longrightarrow \Ext_R^i(Y,k) \longrightarrow
\Ext_R^i(\tau N,k) \stackrel{\delta}{\longrightarrow}$.  Part (i) of
the lemma will follow if the connecting maps $\delta$ are zero.  

And they are: By \cite[lem.\ 4.2]{PJARtq}, the \ARt\ is also an \ARt\
in $\DfR$, the full subcategory of $\D(R)$ consisting of $M$'s with
$\dim_k \H\!M < \infty$.  Moreover, no morphism $\Sigma^{-i}(\tau N)
\longrightarrow {}_{R}k$ is a section, for any such section would 
clearly have to be an isomorphism, but $\Sigma^{-i}(\tau N)$ is in
$\DcR$ and ${}_{R}k$ is not.  The definition of \ARt\ now says that
the composition $\Sigma^{-i-1}N \longrightarrow \Sigma^{-i}(\tau N)
\longrightarrow {}_{R}k$ is zero, where the first arrow is a
suspension of the connecting morphism from the \ARt.  

In other words, the map 
\[
  \Hom_{\Dsmall(R)}(\Sigma^{-i}(\tau N),k)
  \longrightarrow \Hom_{\Dsmall(R)}(\Sigma^{-i-1}N,k)
\]
is zero.  But this map equals
\[
  \Ext_R^i(\tau N,k) 
  \stackrel{\delta}{\longrightarrow} \Ext_R^{i+1}(N,k).
\]

\medskip
\noindent
(ii) Consider again $\tau N \longrightarrow Y \longrightarrow N
\longrightarrow$, the \ARt\ from (i).  From \cite[prop.\
3.5]{HapDerCat} follows that $Y$ is a coproduct of copies of those
indecomposable objects $M$ of $\DcR$ for which there are irreducible
morphisms $M \longrightarrow N$.  That is, $Y$ is a coproduct of
copies of those indecomposable objects $M$ of $\DcR$ for which there are
arrows $[M] \longrightarrow [N]$ in the \ARq\ of $\DcR$.  And by
the definition at the beginning of section \ref{sec:labelling},
the multiplicity of $M$ as a direct summand of $Y$ is
$a_{[M] \rightarrow [N]}$.  This clearly implies
\[
  \AddFunct(Y) =
  \sum_{[M] \rightarrow [N]}a_{[M] \rightarrow [N]}\AddFunct([M]),
\]
and then (ii) follows from (i).

\medskip
\noindent
(iii)  Since both $\DcR$ and $\DcRo$ have \ARt s, it follows from
\cite[cor.\ 5.2 and thm.\ 5.1 and its proof]{PJARtq} that there is an
isomorphism ${}_{R}(\dual\!R) \cong {}_{R}(\Sigma^d R)$ in $\DcR$.  So
as $k$-vector spaces,
\[
  \H(\RHom_R(\dual\!R,k))
  \cong \H(\RHom_R(\Sigma^d R,k))
  \cong \Sigma^{-d}k,
\]
which easily implies $\RHom_R(\dual\!R,k) \cong \Sigma^{-d}k$ in
$\D(R)$.  This again gives (a) in
\begin{align*}
  \RHom_R(\tau N,k) 
  & = \RHom_R(\Sigma^{-1}(\dual\!R \LTensor_R N),k) \\
  & \cong \Sigma \RHom_R(\dual\!R \LTensor_R N,k) \\
  & \cong \Sigma \RHom_R(N,\RHom_R(\dual\!R,k)) \\
  & \stackrel{\rm (a)}{\cong} \Sigma \RHom_R(N,\Sigma^{-d}k) \\
  & \cong \Sigma^{1-d} \RHom_R(N,k),
\end{align*}
and taking cohomology and $k$-dimension shows $\AddFunct(\tau[N]) =
\AddFunct([N])$. 
\end{proof}

The following results leave $\AddFunct$ for a while, but return to it
in corollary \ref{cor:AddFunct_unbounded}.

There is an abelian category $\sC(R)$ whose objects are \DGRlm s and
whose morphisms are homomorphisms of \DGRlm s, that is, homomorphisms
of \grlm s which are compatible with the differentials.  The following
Harada-Sai lemma for $\sC(R)$ is proved just like
\cite[lem.\ 4.14.1]{Bensonbook}.

\begin{Lemma}
\label{lem:composition}
Let $F_0, \ldots, F_{2^p-1}$ be indecomposable objects of $\sC(R)$ with
$\dim_k F_i \leq p$ for each $i$, and let
\[
  \begin{diagram}[labelstyle=\scriptstyle,width=3ex]
    F_{2^p-1} & \rTo^{\varphi_{2^p-1}} & F_{2^p-2} 
    & \rTo^{\varphi_{2^p-2}} & \cdots & \rTo^{\varphi_1} & F_0 \\
  \end{diagram}
\]
be non-isomorphisms in $\sC(R)$.  

Then $\varphi_1 \circ \cdots \circ \varphi_{2^p-1} = 0$.
\end{Lemma}

\begin{Corollary}
\label{cor:composition}
Let $M_0, \ldots, M_{2^p-1}$ be indecomposable objects of $\DcR$ with
$\AddFunct(M_i) \leq \frac{p}{\dim_k R}$ for each $i$, and let 
\[
  \begin{diagram}[labelstyle=\scriptstyle,width=3ex]
    M_{2^p-1} & \rTo^{\mu_{2^p-1}} & M_{2^p-2} 
    & \rTo^{\mu_{2^p-2}} & \cdots & \rTo^{\mu_1} & M_0 \\
  \end{diagram}
\]
be non-isomorphisms in $\DcR$.  

Then $\mu_1 \circ \cdots \circ \mu_{2^p-1} = 0$.
\end{Corollary}

\begin{proof}
Pick minimal semi-free resolutions $F_i
\stackrel{\simeq}{\longrightarrow} M_i$ (minimality means that the
differential $\partial_{F_i}$ takes values inside $R^{\geq 1} \cdot
F_i$, see \cite[lem.\ 3.3]{PJARtq}).  Then the morphisms 
$
\begin{diagram}[labelstyle=\scriptstyle,width=3ex]
  M_{2^p-1} & \rTo^{\mu_{2^p-1}} & M_{2^p-2} 
  & \rTo^{\mu_{2^p-2}} & \cdots & \rTo^{\mu_1} & M_0 \\
\end{diagram}
$
in $\DcR$ are represented by morphisms 
$
\begin{diagram}[labelstyle=\scriptstyle,width=3ex]
  F_{2^p-1} & \rTo^{\varphi_{2^p-1}} & F_{2^p-2} 
  & \rTo^{\varphi_{2^p-2}} & \cdots & \rTo^{\varphi_1} & F_0 \\
\end{diagram}
$ 
in $\sC(R)$.

If $\varphi_i$ were an isomorphism in $\sC(R)$, then $\mu_i$ would be
an isomorphism in $\DcR$, so each $\varphi_i$ is a non-isomorphism in
$\sC(R)$.  Also, since $F_i$ is minimal, it is easy to see that if
$F_i$ decomposed non-trivially in $\sC(R)$, then $M_i$ would decompose
non-trivially in $\DcR$, so each $F_i$ is indecomposable in $\sC(R)$.

Now
\begin{align}
\nonumber
  \Ext_R(M_i,k) & = \H(\RHom_R(M_i,k)) \cong \H(\Hom_R(F_i,k)) \\
\label{equ:Ext}
  & \stackrel{\rm (a)}{\cong} \Hom_R(F_i,k)^{\natural}
  = \Hom_{R^{\natural}}(F_i^{\natural},k^{\natural}).
\end{align} 
Here $\natural$ indicates the functor which forgets differentials, and
therefore sends DGAs to graded algebras and \DGm s to graded modules.
The isomorphism (a) holds because $F_i$ is minimal, whence the
differential of the complex $\Hom_R(F_i,k)$ is zero, so taking
cohomology amounts simply to forgetting the differential.  

Taking $k$-dimensions in equation \eqref{equ:Ext} shows that
$\AddFunct(M_i)$ equals the $k$-dimension of
$\Hom_{R^{\natural}}(F_i^{\natural},k^{\natural})$.  However, since
$F_i$ is semi-free I have
\[
  F_i^{\natural} \cong \bigoplus_j \Sigma^{\sigma_j}(R^{\natural})
\]
for certain integers $\sigma_j$, and the $k$-dimension of
$\Hom_{R^{\natural}}(F_i^{\natural},k^{\natural})$ clearly equals the
number of summands $\Sigma^{\sigma_j}(R^{\natural})$ in
$F_i^{\natural}$.  

All in all, $\AddFunct(M_i)$ equals the number of summands
$\Sigma^{\sigma_j}(R^{\natural})$ in $F_i^{\natural}$ which gives (b)
in
\begin{align*}
  \dim_k F_i & = \dim_k F_i^{\natural}
  \stackrel{\rm (b)}{=} \AddFunct(M_i) \dim_k R^{\natural} \\
  & = \AddFunct(M_i) \dim_k R 
  \leq \frac{p}{\dim_k R} \dim_k R 
  = p.
\end{align*}

But now lemma \ref{lem:composition} gives $\varphi_1 \circ \cdots
\circ \varphi_{2^p - 1} = 0$ which clearly implies $\mu_i \circ
\cdots \circ \mu_{2^p - 1} = 0$.  
\end{proof}

\begin{Lemma}
\label{lem:non_zero_composition}
Suppose that $\DcR$ has \ARt s and that ${}_{R}k$ is not in $\DcR$.

Given an indecomposable object $M_0$ of $\DcR$ and an integer $q \geq
0$, there exist indecomposable objects and irreducible morphisms
\[
  \begin{diagram}[labelstyle=\scriptstyle,width=3ex]
    M_q & \rTo^{\mu_q} & M_{q-1} & \rTo^{\mu_{q-1}}
    & \cdots & \rTo^{\mu_1} & M_0 \\
  \end{diagram}
\]
in $\DcR$ with $\mu_1 \circ \cdots \circ \mu_q \not= 0$.
\end{Lemma}

\begin{proof}
If $M$ is a non-zero object in $\DfR$, the full subcategory of $\D(R)$
consisting of $M$'s with $\dim_k \H\!M < \infty$, then the $k$-dual
$\dual\!M = \Hom_k(M,k)$ is non-zero in $\DfRo$, so the minimal
semi-free resolution $F$ of $\dual\!M$ is non-trivial.  Hence
$\RHom_{R^{\opp}}(\dual\!M,k) \cong \Hom_{R^{\opp}}(F,k)$ has non-zero
cohomology, so by dualization the same holds for
\[
  \RHom_R(\dual\!k,\dual\!\dual\!M) \cong \RHom_R(k,M).  
\]
So there exists an $i$ so that $\H^{-i}(\RHom_R(k,M)) \cong
\Hom_{\Dsmall(R)}(\Sigma^i k,M)$ is non-zero;
that is,
\begin{eqnarray}
\nonumber
  & \mbox{For $M$ non-zero in $\DfR$, there exists} & \\
\label{equ:non_zero_morphism}
  & \mbox{a non-zero morphism 
          ${}_{R}(\Sigma^i k) \longrightarrow M$ in $\DfR$.} &
\end{eqnarray}

Moreover, since $M$ is non-zero, each retraction ${}_{R}(\Sigma^i k)
\longrightarrow M$ in $\DfR$ is clearly an isomorphism.  If $M$ is in
$\DcR$, then there is no such isomorphism because ${}_{R}(\Sigma^i k)$
is not in $\DcR$, so 
\begin{eqnarray}
\nonumber
  & \mbox{For $M$ non-zero in $\DcR$, no} & \\
\label{equ:non_retraction}
  & \mbox{morphism ${}_{R}(\Sigma^i k) \longrightarrow M$ 
          is a retraction in $\DfR$.} &
\end{eqnarray}

Now for the proof proper.  In fact, I shall prove slightly more than
claimed, namely, there exists
\[
  \begin{diagram}[labelstyle=\scriptstyle,width=3ex]
    {}_{R}(\Sigma^i k) & \rTo^{\kappa_q}
    & M_q & \rTo^{\mu_q} & M_{q-1} & \rTo^{\mu_{q-1}}
    & \cdots & \rTo^{\mu_1} & M_0 \\
  \end{diagram}
\]
with the $M_i$ indecomposable and the $\mu_i$ irreducible in $\DcR$,
with $\kappa_q$ not a retraction in $\DfR$, and with $\mu_1 \circ
\cdots \circ \mu_q \circ \kappa_q \not= 0$.  I do so by induction on
$q$. 

For $q = 0$, existence of $\kappa_0$ holds by equations
\eqref{equ:non_zero_morphism} and \eqref{equ:non_retraction}. 

For $q \geq 1$, the induction gives data
\[
  \begin{diagram}[labelstyle=\scriptstyle,width=3ex,midshaft]
    {}_{R}(\Sigma^i k) & \rTo^{\kappa_{q-1}}
    & M_{q-1} & \rTo^{\mu_{q-1}} &
    & \cdots & \rTo^{\mu_1} & M_0. \\
  \end{diagram}
\]
Let
$
\begin{diagram}[labelstyle=\scriptstyle,width=3ex]
\tau M_{q-1} & \rTo & X_q & \rTo^{\mu_q^{\prime}} & M_{q-1} & \rTo \\
\end{diagram}
$ 
be an \ARt\ in $\DcR$.  By \cite[lem.\ 4.2]{PJARtq} this is also an
\ARt\ in $\DfR$, so since
$
\begin{diagram}[labelstyle=\scriptstyle,width=3ex]
{}_{R}(\Sigma^i k) & \rTo^{\kappa_{q-1}} & M_{q-1} \\
\end{diagram}
$ 
is not a retraction in $\DfR$, it factors through $X_q$ as 
$
\begin{diagram}[labelstyle=\scriptstyle,width=3ex]
{}_{R}(\Sigma^i k) & \rTo^{\kappa_q^{\prime}} & X_q
& \rTo^{\mu_q^{\prime}} & M_{q-1}, \\
\end{diagram}
$
and so $\mu_1 \circ \cdots \circ \mu_{q-1} \circ \mu_q^{\prime} \circ
\kappa_q^{\prime} = \mu_1 \circ \cdots \circ \mu_{q-1} \circ
\kappa_{q-1}$ is non-zero.

Splitting $X_q$ into indecomposable summands and $\kappa_q^{\prime}$
and $\mu_q^{\prime}$ into components gives that there exists 
$
\begin{diagram}[labelstyle=\scriptstyle,width=3ex]
{}_{R}(\Sigma^i k) & \rTo^{\kappa_q} & M_q & \rTo^{\mu_q} & M_{q-1} \\
\end{diagram}
$ 
with $\mu_1 \circ \cdots \circ \mu_{q-1} \circ \mu_q \circ \kappa_q
\not= 0$.  

Here $M_q$ is an indecomposable summand of $X_q$ so is in $\DcR$, so
equation \eqref{equ:non_retraction} gives that $\kappa_q$ is not a
retraction in $\DfR$.  And $\mu_q$ is a component of $\mu_q^{\prime}$, so
$\mu_q$ is irreducible in $\DcR$ by \cite[prop.\ 3.5]{HapDerCat}.
\end{proof}

I now return to the function $\AddFunct$.

\begin{Corollary}
\label{cor:AddFunct_unbounded}
Suppose that $\DcR$ and $\DcRo$ have \ARt s and that ${}_{R}k$ is not
in $\DcR$.  Let $\Component$ be a component of the \ARq\ of $\DcR$.
Then $\AddFunct$ is unbounded on $\Component$.
\end{Corollary}

\begin{proof}
Let $M_0$ be an indecomposable object of $\DcR$ for which $[M_0]$ is a
vertex of $C$.  Lemma \ref{lem:non_zero_composition} says that there
exist arbitrarily long sequences of indecomposable objects and
irreducible morphisms
\[
  \begin{diagram}[labelstyle=\scriptstyle,width=3ex]
    M_q & \rTo^{\mu_q} & M_{q-1} & \rTo^{\mu_{q-1}}
    & \cdots & \rTo^{\mu_1} & M_0 \\
  \end{diagram}
\]
in $\DcR$ with $\mu_1 \circ \cdots \circ \mu_q \not= 0$.  Each $[M_i]$
is clearly a vertex of $\Component$.  So corollary
\ref{cor:composition} implies that $\AddFunct$ cannot be bounded on
$\Component$.
\end{proof}

\section{Structure of the \ARq}
\label{sec:shape}

This section combines the material of previous sections to prove the
main result of this paper, theorem \ref{thm:topology}.

By lemma \ref{lem:Riedtmann}, if $\DcR$ and $\DcRo$ have \ARt s, then
each component $\Component$ of the \ARq\ of $\DcR$ is of the form $\BZ
\Tree/\Group$, where $\Tree$ is a directed tree and $\Group \subseteq
\Aut(\BZ \Tree)$ is an admissible group of automorphisms.

Recall from \cite[sec.\ 4.15]{Bensonbook} that the vertices of $\BZ
\Tree$ are the pairs $(p,t)$ where $p$ is an integer and $t$ is a
vertex of $\Tree$, and that the vertices of $\BZ \Tree/\Group$ are the
orbits $\Group(p,t)$ under the group $\Group$.  By identifying $\BZ
\Tree/\Group$ with $\Component$, I shall also view the $\Group(p,t)$'s
as vertices of $\Component$, so the following makes sense.

First, the start of section \ref{sec:labelling} constructs a labelling
of the \ARq, and so in particular a labelling of $\Component$.  Now,
if there is an arrow $t \longrightarrow u$ in $\Tree$, then there is
an arrow $\Group(0,t) \longrightarrow \Group(0,u)$ in $\BZ
\Tree/\Group$, that is, an arrow $\Group(0,t) \longrightarrow
\Group(0,u)$ in $\Component$.  Hence the labelling of $\Component$
induces a labelling of $\Tree$ by
\begin{equation}
\label{equ:Tree_labelling}
  (a_{t \rightarrow u},b_{t \rightarrow u}) =
  (a_{\Group(0,t) \rightarrow \Group(0,u)},
   b_{\Group(0,t) \rightarrow \Group(0,u)}).
\end{equation}

Secondly, the start of section \ref{sec:AddFunct} defines a function
$\AddFunct$ on the vertices of $\Component$.  That is, $\AddFunct$ is
defined on the $\Group(p,t)$'s, and so induces a function
$\AddFunctTree$ on $\Tree$ by
\begin{equation}
\label{equ:AddFunctTree}
  \AddFunctTree(t) = \AddFunct(\Group(0,t)).
\end{equation}

\begin{Lemma}
\label{lem:induction}
Suppose that $\DcR$ and $\DcRo$ have \ARt s.
\begin{enumerate}

  \item  The function $\AddFunctTree$ is additive with respect to the
         labelling of $\Tree$ given by equation
         \eqref{equ:Tree_labelling}.

  \item  The function $\AddFunctTree$ is unbounded on $\Tree$.

\end{enumerate}
\end{Lemma}

To explain part (i), let me recall from \cite[sec.\ 4.5]{Bensonbook}
that when $\Tree$ is a labelled directed tree, a function
$\AddFunctTree$ on the vertices of $\Tree$ is called additive if it
satisfies
\begin{equation}
\label{equ:f_additive}
  2\AddFunctTree(t) 
  - \sum_{s \rightarrow t}a_{s \rightarrow t}\AddFunctTree(s)
  - \sum_{t \rightarrow u}b_{t \rightarrow u}\AddFunctTree(u) = 0
\end{equation}
for each vertex $t$. 

\begin{proof}
[Proof of Lemma \ref{lem:induction}]
(i)  To show that the left hand side of equation \eqref{equ:f_additive} 
is zero, let me rewrite it as
\begin{align}
\nonumber
  2\AddFunct(\Group(0,t))
  & - \sum_{s \rightarrow t}
        a_{\Group(0,s) \rightarrow \Group(0,t)}\AddFunct(\Group(0,s)) \\
\label{equ:form_i}
  & - \sum_{t \rightarrow u}
        b_{\Group(0,t) \rightarrow \Group(0,u)}\AddFunct(\Group(0,u)),
\end{align}
using equations \eqref{equ:Tree_labelling} and \eqref{equ:AddFunctTree}.

Recall from \cite[sec.\ 4.15]{Bensonbook} that 
the translation of $\BZ \Tree/\Group$ is
given by $\tau(\Group(p,t)) =
\Group(p+1,t)$.  Identifying $\BZ \Tree/\Group$ with $\Component$,
lem\-ma \ref{lem:symmetry}(i) gives 
$b_{\Group(0,t) \longrightarrow \Group(0,u)} 
= a_{\tau(\Group(0,u)) \longrightarrow \Group(0,t)} 
= a_{\Group(1,u) \longrightarrow \Group(0,t)}$,
and lemma \ref{lem:AddFunct}(iii) gives
$\AddFunct(\Group(0,v)) 
= \AddFunct(\tau(\Group(0,v))) 
= \AddFunct(\Group(1,v))$.
Hen\-ce \eqref{equ:form_i} is 
\begin{align}
\nonumber
  \AddFunct(\tau \Group(0,t)) + \AddFunct(\Group(0,t))
  & - \sum_{s \rightarrow t}
        a_{\Group(0,s) \rightarrow \Group(0,t)}\AddFunct(\Group(0,s)) \\
\label{equ:form_ii}
  & - \sum_{t \rightarrow u}
        a_{\Group(1,u) \rightarrow \Group(0,t)}\AddFunct(\Group(1,u)).
\end{align}

To rewrite further, recall also from \cite[sec.\ 4.15]{Bensonbook}
that the arrows in $\BZ \Tree$ which end in $(0,t)$ are obtained as
follows: There is an arrow $(0,s) \longrightarrow (0,t)$ in $\BZ
\Tree$ for each arrow $s \longrightarrow t$ in $\Tree$, and there is
an arrow $(1,u) \longrightarrow (0,t)$ in $\BZ \Tree$ for each arrow
$t \longrightarrow u$ in $\Tree$.

Moreover, the canonical map $\BZ \Tree \longrightarrow \BZ
\Tree/\Group$ is a so-called covering, see \cite[p.\ 156]{Bensonbook},
so the map which sends the arrow $m \longrightarrow (0,t)$ to the
arrow $\Group m \longrightarrow \Group(0,t)$ is a bijection between
the arrows in $\BZ \Tree$ which end in $(0,t)$ and the arrows in $\BZ
\Tree/\Group$ which end in $\Group(0,t)$.  All this implies that taken
together, the two sums in \eqref{equ:form_ii} contain exactly one
summand for each arrow in $\BZ \Tree/\Group$ which ends in
$\Group(0,t)$, so \eqref{equ:form_ii} is
\[
  \AddFunct(\tau \Group(0,t)) + \AddFunct(\Group(0,t))
  - \sum_{m \rightarrow \Group(0,t)}
      a_{m \rightarrow \Group(0,t)}\AddFunct(m).
\]
Identifying $\BZ \Tree/\Group$ with $\Component$, this is zero by
lemma \ref{lem:AddFunct}(ii).

\medskip
\noindent
(ii)  I have 
\[
  \AddFunct(\Group(p,t)) = \AddFunct(\tau^p(\Group(0,t)))
  \stackrel{\rm (a)}{=} \AddFunct(\Group(0,t)) = \AddFunctTree(t),
\]
where (a) is by lemma \ref{lem:AddFunct}(iii).  So if $\AddFunctTree$
were bounded on $\Tree$, then $\AddFunct$ would be bounded on $C$
which it is not by corollary \ref{cor:AddFunct_unbounded}.
\end{proof}

The following is my abstract main result, which has theorem
\ref{thm:topology} as an easy consequence.

I emphasize that as above, $\Component$ is a component of the \ARq\
viewed as a stable translation quiver, so $\Component$ is closed under
the \ARtr\ $\tau$ and its inverse, and is a connected stable
translation quiver with translation the restriction of $\tau$.

\begin{Theorem}
\label{thm:main}
Suppose that $\DcR$ and $\DcRo$ have \ARt s and that ${}_{R}k$ is not
in $\DcR$.  Let $\Component$ be a component of the \ARq\ of $\DcR$.

Then $\Component$ is isomorphic to $\BZ A_{\infty}$ as a stable
translation quiver.
\end{Theorem}

\begin{proof}
Lemma \ref{lem:Riedtmann} says that $\Component$ is of the form $\BZ
\Tree/\Group$. 

Lemma \ref{lem:induction} gives an additive function $\AddFunctTree$
on $\Tree$, and shows that $\AddFunctTree$ is unbounded.
Consequently, the underlying graph of $\Tree$ must be $A_{\infty}$ by
\cite[thm.\ 4.5.8(iv)]{Bensonbook}.  So $\Component$ is $\BZ
A_{\infty}/\Group$.  (Note that when the underlying graph of $\Tree$
is $A_{\infty}$, the quiver $\BZ \Tree$ is determined up to
isomorphism, independently of the directions of the arrows in
$\Tree$.  So it makes sense to denote $\BZ \Tree$ by $\BZ
A_{\infty}$.) 

Now, $\Group$ must send vertices at the end of $\BZ A_{\infty}$ to
other such vertices.  If $\Group$ acted non-trivially on the vertices
at the end of $\BZ A_{\infty}$, then there would exist a $g$ in
$\Group$ and a vertex $m$ at the end of $\BZ A_{\infty}$ so that $gm$
was a different vertex at the end of $\BZ A_{\infty}$.  As the other
vertices at the end of $\BZ A_{\infty}$ have the form $\tau^p m$ for
$p \not= 0$, this would mean $gm = \tau^p m$ for some $p \not= 0$.
Then $m$ and $\tau^p m$ would get identified in $\BZ
A_{\infty}/\Group$, so $\Group m$ would be a fixed point of $\tau^p$
in $\BZ A_{\infty}/\Group$, that is, a fixed point of $\tau^p$ in
$\Component$ contradicting lemma \ref{lem:tau_automorphism}(iii).

So $\Group$ acts trivially on the vertices at the end of $\BZ
A_{\infty}$, and it is easy to see that this forces $\Group$ to act
trivially on all of $\BZ A_{\infty}$.  So $\BZ A_{\infty}/\Group$ is
just $\BZ A_{\infty}$, and the theorem follows.
\end{proof}

\begin{proof}
[Proof of Theorem \ref{thm:topology}]
This proof uses a bit of rational homotopy theory for which my source
is \cite{FHTbook}.

By \cite[exam.\ 6, p.\ 146]{FHTbook} I have that $\CXk$ is equivalent
by a series of quasi-isomorphisms to a DGA, $R$, which satisfies setup
\ref{set:blanket}.  Hence the various derived categories of $\CXk$ are
equivalent to those of $R$ by \cite[thm.\ III.4.2]{KrizMayAst}, so
theorem \ref{thm:main} can be applied to $\CXk$.

To get the desired conclusion, that each component $\Component$ of the
\ARq\ of $\DcXk$ is isomorphic to $\BZ A_{\infty}$, I must show that
the premises of theorem \ref{thm:main} hold for $\CXk$.  So I must
show that $\Dc(\CXk)$ and $\Dc(\CXk^{\opp})$ have \ARt s, and that
${}_{\CXk}k$ is not in $\Dc(\CXk)$.  The first two claims hold by the
main result of \cite{PJARtq}, theorem 6.3, because $\Space$ has
Poincar\'e duality over $k$.  The last claim can be seen as follows.

By assumption, $\Space$ has Poincar\'e duality of dimension $d \geq 2$
over $k$, so I have $\H^{\geq 2}(\Space;k) \not= 0$.  Hence the
minimal Sullivan model $\Lambda V$ of $\CXk$ has $V^{\geq 2} \not= 0$;
see \cite[chp.\ 12]{FHTbook}.  But $V$ is
$\Hom_{\BZ}(\pi_{\ast}\Space,k)$ by \cite[chp.\ 15]{FHTbook}, where
$\pi_{\ast}\Space$ is the sequence of homotopy groups of $\Space$, so
as $k$ has characteristic zero, $\pi_{\geq 2}\Space$ cannot be
torsion, so $\BQ \otimes_{\BZ}
\pi_{\geq 2}\Space \not= 0$.  Since \cite[p.\ 434]{FHTbook} gives
$\dim_{\BQ}(\BQ \otimes_{\BZ} \pi_{\rm odd}\Space) -
\dim_{\BQ} (\BQ \otimes_{\BZ} \pi_{\rm even}\Space) \geq 0$, there
follows $\BQ \otimes_{\BZ} \pi_{\rm odd}\Space \not= 0$.  This implies
that the denominator is non-trivial in the formula
\cite[(33.7)]{FHTbook} from which follows $\dim_{\BQ} \H^{\ast}(\Omega
\Space;\BQ) = \infty$, where $\Omega \Space$ is the Moore loop space
of $\Space$.  So also
\begin{equation}
\label{equ:H_Omega_X_infinite}
  \dim_k \H^{\ast}(\Omega \Space;k) = \infty.
\end{equation}

Now consider the Moore path space fibration $\Omega \Space
\longrightarrow P\Space \longrightarrow \Space$ from \cite[exam.\ 1,
p.\ 29]{FHTbook}.
The Moore path space $P\Space$ is contractible, and this implies that the
\DGCXklm\ $\C^{\ast}(P\Space;k)$ is quasi-isomorphic to
${}_{\CXk}k$.  Inserting the fibration
into the Eilenberg-Moore theorem \cite[thm.\
7.5]{FHTbook} therefore gives a quasi-isomorphism
$\C^{\ast}(\Omega \Space;k) \simeq k \LTensor_{\CXk} k$,
and taking cohomology shows 
\[
  \H^{\ast}(\Omega \Space;k) 
  \cong \Tor_{-\ast}^{\CXk}(k,k).  
\]
Com\-bi\-ning with equation \eqref{equ:H_Omega_X_infinite}
gives 
\[
  \dim_k \Tor^{\CXk}(k,k) = \infty.
\]

So ${}_{\CXk}k$ cannot be finitely built from ${}_{\CXk}\CXk$ whence
${}_{\CXk}k$ is not in $\Dc(\CXk)$ by \cite[lem.\ 3.2]{PJARtq}.
\end{proof}

\end{document}